\documentclass{amsart}

\usepackage{url}
\usepackage{amssymb, amsmath, amscd, amsthm}
\usepackage{graphicx}
\usepackage{psfrag}
\usepackage[all]{xy}
\usepackage{rotating} 

 \newtheorem{teo}{\bf Theorem}[section]
 \newtheorem{thm}[teo]{\bf Theorem}
 \newtheorem{lem}[teo]{\bf Lemma}
 
 \newtheorem{prop}[teo]{\bf Proposition}

 \newtheorem{defi}[teo]{\bf Definition}
 
 \newtheorem{rem}[teo]{\bf Remark}

 \theoremstyle{definition}
 \newtheorem{exa}{\bf Example}

\newcommand{\ba}{\begin{array}}
\newcommand{\ea}{\end{array}}

\newcommand{\R}{\mathbb{R}}
\newcommand{\Q}{\mathbb{Q}}

\newcommand{\f}{f_*}

\newcommand{\p}{\varphi}
\newcommand{\s}{\psi}

\newcommand{\eps}{\varepsilon}

\begin{document}

\title[Optimal homeomorphisms between closed curves]{Optimal homeomorphisms between closed curves}
\author[A. Cerri]{Andrea Cerri}
\address{ARCES, Universit\`a di Bologna,
via Toffano $2/2$, I-$40135$ Bologna, Italia\newline Dipartimento
di Matematica, Universit\`a di Bologna, P.zza di Porta S. Donato
5, I-$40126$ Bologna, Italia}
\email{\{cerri,difabio\}@dm.unibo.it}

\author[B. Di Fabio]{Barbara Di Fabio}


\date{\today} 



\begin{abstract}
The concept of  natural pseudo-distance has proven to be a
powerful tool for measuring the dissimilarity between topological
spaces endowed with continuous real-valued functions. Roughly
speaking, the natural pseudo-distance is defined as the infimum of
the change of the functions' values, when moving from one space to
the other through homeomorphisms, if possible. In this paper, we
prove the first available result about the existence of optimal
homeomorphisms between closed curves, i.e. inducing a change of the
function that equals the natural pseudo-distance.
\end{abstract}

\maketitle

\section*{Introduction}

Formalizing the concept of shape for topological spaces and
manifolds, as well as providing an efficient comparison of shapes,
has been a widely researched topic in the last decade. As such, a
class of methods has been developed with the purpose of performing
a topological exploration of the shape, according to some
quantitative geometric properties provided by a real function
chosen to extract shape features
\cite{BiDeFaFrGiLaPaSp08,Ca09,Gh08,KaMiMr04,Zo05}.

In this context, Size Theory was introduced at the beginning of
the 1990s \cite{Fr90,Fr91,FrLa99}, 
supported by the adoption of a suitable mathematical tool: the
{\em natural pseudo-distance} \cite{DoFr94,DoFr04,DoFr07}.

In the formalism of Size Theory, a shape is modelled as a pair
$(X,\p)$, where $X$ is a topological space and $\p:X\to\R$ is a
continuous function \cite{BiDeFaFrGiLaPaSp08,FrLa99}. Such a pair
is called a {\em size pair} and $\p$ is called a {\em measuring
function}. The role of $\p$ is to take into account only the
properties considered relevant for the shape comparison problem at
hand, while disregarding the irrelevant ones, as well as to impose
the desired invariance properties.

The natural pseudo-distance is a measure of the dissimilarity
between two size pairs $(X,\p)$, $(Y,\s)$. Roughly speaking, it is
defined as the infimum of the variation of the values of $\p$ and
$\s$, when we move from $X$ to $Y$ through homeomorphisms, if
possible (see Definition \ref{pseudod}). Therefore, two objects
have the same shape if they share the same shape properties,
expressed by the measuring functions' values, that is, their
natural pseudo-distance vanishes.

Earlier results about the natural pseudo-distance can be divided
into two classes. One class provides constraints on the possible
values taken by the natural pseudo-distance between two size pairs
$(X,\p)$, $(Y,\s)$. For example, if the considered topological
spaces $X$ and $Y$ are smooth closed manifolds and the measuring
functions are also smooth, then the natural pseudo-distance is an
integer sub-multiple of the Euclidean distance between two
suitable critical values of the measuring functions \cite{DoFr04}.
In particular, this integer can only be either 1 or 2 in the case
of curves \cite{DoFr94}, while it cannot be greater than 3 in the
case of surfaces \cite{DoFr07}.

The other class of results furnishes lower bounds for the natural
pseudo-distance \cite{DoFr04bis,FrMu99}. In particular it is
possible to estimate the natural pseudo-distance by using the
concept of \emph{size function} \cite{dAFrLa,DoFr04bis}. Size
functions are shape descriptors able to reduce the comparison of
shapes to the comparison of certain countable sets of points in
the real plane \cite{DoFrLa99,FrLa97,FrLa01}. This reduction
allows us to study the space of all homeomorphisms between the
considered topological spaces, without actually computing them.

The research on size functions has led to a formal setting, which
has turned out to be useful, not only from a theoretical point of
view, but also on the applicative side (see, e.g.,
\cite{BiGiSpFa08,CeFeGi06,DiFrPa04,StBr05,UrVe97}).

Besides being a useful theoretical tool for applications in shape
comparison, the natural pseudo-distance is challenging from the
mathematical point of view, and several questions about its
properties need for further investigation. One among them consists
in establishing the hypotheses ensuring the existence of {\em
optimal homeomorphisms} between size pairs, i.e. homeomorphisms
realizing the natural pseudo-distance. It is possible to show
that, in general, such homeomorphisms do not exist (see, e.g.,
Section \ref{examples}).

In this paper, we provide the first available result about the
existence of optimal homeomorphisms. To be more precise, we prove
that, under appropriate conditions, it is always possible to
construct a homeomorphism between two closed curves (i.e. compact
and without boundary $1$-manifolds), satisfying the property of
optimality (Theorem \ref{mainthm}). This result can be seen as a
necessary first step towards the study of this problem in a more
general setting, e.g. when manifolds of arbitrary dimensions are
involved.

The subject of our work fits in the current mathematical research
and interest in simple closed curves, motivated by problems
concerning shape comparison in Computer Vision (cf. e.g.,
\cite{MiMu06, MiMu07}).

The paper is divided into three sections. Section \ref{background}
deals with some of the standard facts on the comparison of size
pairs via the natural pseudo-distance. In particular, the
definition of the natural pseudo-distance $\delta$ and its main
properties are given, focusing on the concept of $d$-approximating
sequence. Section \ref{examples} is devoted to the description of
some simple and meaningful examples showing that none of the
conditions we require in stating our main result can be dropped.
In Section \ref{optimal} we prove our main result concerning the
existence and the construction of an optimal homeomorphism between
two smooth closed curves endowed with Morse measuring functions
(Theorem \ref{mainthm}).

\section{Preliminaries}\label{background}

In Size Theory, a {\em size pair} is a pair $(X,\p)$, where $X$ is
a non-empty, compact, locally connected Hausdorff space and
$\p:X\to\R$ is a continuous function called a {\em measuring
function}. Let $Size$ be the collection of all the size pairs, and
let $(X,\p),(Y,\s)$ be two size pairs. We denote by $H(X,Y)$ the
set of all homeomorphisms from $X$ to $Y$.
\begin{defi}\label{sizemeasure}
If $H(X,Y)\neq \emptyset$, the function $\Theta: H(X,Y)
\rightarrow \R$ given by
$$\Theta(f)=\underset{x\in X}\max|\p(x)-\s(f(x))|$$
is called the \emph{natural size measure} with respect to the
measuring functions $\p$ and $\s$.
\end{defi}

Roughly speaking, $\Theta(f)$ measures how much $f$ changes the
values taken by the measuring functions, at corresponding points.

\begin{defi}\label{pseudod}
We shall call \emph{natural pseudo-distance} the pseudo-distance
$\delta: Size\times Size\rightarrow\mathbb{R}\cup\{+\infty\}$
defined as
$$
\delta\left((X,\p),(Y,\s)\right)=\left\{
\begin{array}{ll}
\underset{f\in H(X, Y)}\inf\Theta(f), & \mbox{if} \,\, H(X, Y)\neq
\emptyset \\ +\infty, & \mbox{otherwise.}
\end{array}\right.
$$
\end{defi}

Note that $\delta$ is not a distance, since two different size
pairs $(X,\p),(Y,\s)$ can have a vanishing pseudo-distance. In
that case, $X$ and $Y$ are only sharing the same shape properties
with respect to the chosen functions $\p$ and $\s$, respectively.
Moreover, we observe that the infimum of $\Theta(f)$ for $f$
varying in $H(X,Y)$ is not always attained. When it is, we shall
say that each homeomorphism $f\in H(X,Y)$ with
$\Theta(f)=\delta\left((X,\p),(Y,\s)\right)$ is an {\em optimal
homeomorphism}. On the other hand, Definition \ref{pseudod}
implies that, if $H(X,Y)\neq\emptyset$, we can always find a
sequence $(f_k)$ of homeomorphisms from $X$ to $Y$, such that
$\underset{k\to\infty}\lim\Theta(f_k)=\delta\left((X,\p),(Y,\s)\right)$.

\begin{defi}\label{DeSeApp}
Let $(X,\p), (Y,\s)$ be two size pairs, with $X,Y$ homeomorphic
and $\delta\left((X,\p),(Y,\s)\right)=d$. Every sequence $(f_k)$
of homeomorphisms $f_k:X\to Y$ such that
$\underset{k\to\infty}\lim\Theta(f_k)=d$ is said to be a
\emph{$d$-approximating sequence} from $(X,\p)$ to $(Y,\s)$.
\end{defi}

\begin{rem}\label{RemDesApp}
We observe that $(f_k)$ is a $d$-approximating sequence from
$(X,\p)$ to $(Y,\s)$ if and only if $(f_k^{-1})$ is a
$d$-approximating sequence from $(Y,\s)$ to $(X,\p)$.
\end{rem}

The main goal of this paper is to show that an optimal
homeomorphism exists between two size pairs $(X,\p)$ and $(Y,\s)$,
under the following conditions:
\begin{enumerate}\label{Hypotheses}
    \item[$(a)$] $(X,\p)$ and $(Y,\s)$ have vanishing natural
    pseudo-distance, i.e. it holds that $\delta((X,\p),(Y,\s))=0$;
    \item[$(b)$] $X$ and $Y$ are two curves of class $C^2$;
    \item[$(c)$] $\p$ and $\s$ are Morse (i.e., smooth and having
    invertible Hessian at each critical point) measuring
    functions.
\end{enumerate}

This result will be formally given and proved later (Theorem
\ref{mainthm}), in the case of closed curves. However, we remark
that the hypothesis of closed curves will be assumed only for the
sake of simplicity. Indeed, it can be weakened to compact
$1$-manifolds having non-empty boundary, without much affecting
the following reasonings.

\begin{rem}
The reader may wonder why we are defining the natural
pseudo-distance $\delta$ in terms of homeomorphisms instead of
diffeomorphisms, since the above assumption $(c)$ requires that
the measuring functions are Morse. The answer is that, under
condition $(b)$, Definition \ref{pseudod} is invariant with
respect to such a choice. Indeed, it is well known that each
homeomorphism between compact differentiable $1$-manifolds can be
approximated arbitrarily well by diffeomorphisms. On the other
hand, dealing with homeomorphisms allows us to slim down examples
and proves from useless technical steps.
\end{rem}

\section{Meaningful examples}\label{examples}

We provide here three meaningful examples showing that the
assumptions $(a)$, $(b)$, $(c)$ introduced in Section
\ref{background} are the less restrictive we can consider in order
to ensure the existence of an optimal homeomorphism between two
size pairs $(X,\p)$ and $(Y,\s)$. Indeed, if one among them is
dropped, then Theorem \ref{mainthm} does not hold.

First of all, let us observe that for every size pair $(Z,\omega)$
with $Z$ a closed curve of class $C^k$, an embedding $h:Z\to\R^3$
of class $C^k$ exists such that $z(p)=\omega(h^{-1}(p))$ for each
point $p=(x(p),y(p),z(p))\in h(Z)$. Moreover, if $\omega$ is
Morse, we can assume that $z$ is Morse on $h(Z)$, too. In other
words, there is no lack of generality in assuming that the
measuring function associated with $Z$ is obtained by restriction
of the $z$-coordinate in $\R^3$.

Accordingly, in the examples and figures we describe here, we
shall always assume that the spaces $X$ and $Y$ are endowed with
the $z$-coordinate function, and use the symbol $z$ to denote both
$z_{|_X}$ and $z_{|_Y}$.

\begin{exa}[Hypothesis $(a)$ fails]
We report an example introduced in \cite{DoFr04}. It shows that if two
size pairs satisfy hypotheses $(b)$ and $(c)$, but have
non-vanishing natural pseudo-distance, then an optimal
homeomorphism does not always exist.

\begin{figure}[htbp]
\begin{center}
\psfrag{A}{$x_A$} \psfrag{B}{$x_B$}
\psfrag{C}{$x_C$}
\psfrag{D}{$x_D^{\,\eps}$}
\psfrag{E}{$x_E^{\,\eps}$} \psfrag{F}{$y_E^{\,\eps}$}
\psfrag{G}{$y_C$} \psfrag{H}{$y_D^{\,\eps}$} \psfrag{g}{$g_{\eps}$}
\psfrag{z}{$z$} \psfrag{e}{$2\eps$} \psfrag{X}{$X$} \psfrag{Y}{$Y$}
\psfrag{M}{$\max z$} \psfrag{m}{$\min z$}
\includegraphics[height=6cm]{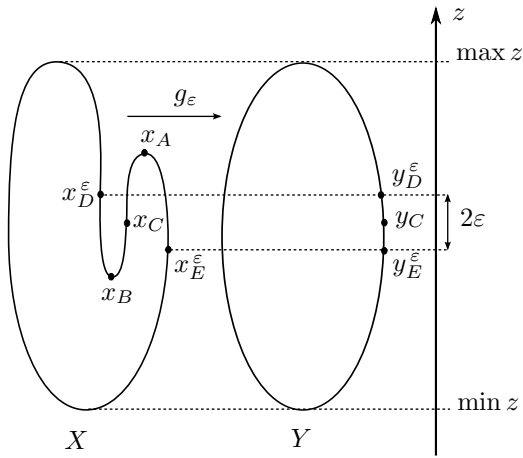}
\caption{\footnotesize{An example of two closed curves $X$, $Y$
endowed with the Morse function $z$. No optimal homeomorphism
exist between $(X, z)$, $(Y, z)$ because their natural
pseudo-distance is non-zero.}}\label{exno(a)}
\end{center}
\end{figure}

Let us consider the two size pairs $(X,z)$, $(Y,z)$ depicted in
Figure $1$, where $X$ and $Y$ are smooth closed curves in $\R^3$,
embedded in the real plane. The functions $z_{|_X}$ and $z_{|_Y}$ are Morse.

As can be seen in Figure \ref{exno(a)}, the points $x_A^{}$,
$x_B^{}\in X$ are critical points of the function $z$ and
$z(x_C^{})=\frac{1}{2}(z(x_A^{})+z(x_B^{}))=z(y_C^{})$. In
\cite{DoFr94} it has been proved that the natural pseudo-distance
between homeomorphic smooth closed curves, endowed with Morse
measuring functions, is always obtainable in terms of some
critical values of the measuring functions. Actually, in this
example it is possible to show that the natural pseudo-distance
between $(X,z)$ and $(Y,z)$ takes the value
$d=\frac{1}{2}(z(x_A^{})-z(x_B^{}))$. On the other hand, it will
also be proved that no optimal homeomorphism exists. Indeed, we
can construct a sequence of homeomorphisms $(f_k)$, such that
$\underset{k\to\infty}\lim\Theta(f_k)=\frac{1}{2}(z(x_A^{})-z(x_B^{}))$,
and show that $\Theta(f)>\frac{1}{2}(z(x_A^{})-z(x_B^{}))$ for
every homeomorphism $f\in H(X,Y)$. The first step consists in
proving that, for every $\varepsilon>0$, a homeomorphism
$g_{\varepsilon}:X\to Y$ exists, such that
$\Theta(g_{\varepsilon})\leq\frac{1}{2}(z(x_A^{})-z(x_B^{}))+2\varepsilon$.
Accordingly, consider the points $x_D^{\,\varepsilon}$,
$x_E^{\,\varepsilon}$, $y_D^{\,\varepsilon}$ and
$y_E^{\,\varepsilon}$ in Figure \ref{exno(a)}, verifying
$z(x_D^{\,\varepsilon})=z(y_D^{\,\varepsilon})=z(x_C^{})+\varepsilon$
and
$z(x_E^{\,\varepsilon})=z(y_E^{\,\varepsilon})=z(x_C^{})-\varepsilon$.
Choose a homeomorphism $g_{\varepsilon}$, taking the arc
$x_D^{\,\varepsilon}\,x_C^{}x_E^{\,\varepsilon}$ to the arc
$y_D^{\,\varepsilon}\,y_C^{}y_E^{\,\varepsilon}$ in such a way
that $g_{\varepsilon}(x_D^{\,\varepsilon})=y_D^{\,\varepsilon}$
and $g_{\varepsilon}(x_E^{\,\varepsilon})=y_E^{\,\varepsilon}$.
Outside the arc $x_D^{\,\varepsilon}x_C^{}x_E^{\,\varepsilon}$ in
$X$ define $g_{\varepsilon}$ by mapping, in the unique possible
way, every point $x$ to a point $g_{\varepsilon}(x)$ satisfying
$z(x)=z(g_{\varepsilon}(x))$. For every
$k\in\mathbb{N}\setminus\{0\}$ set $f_k=g_{\frac{1}{k}}$. It can
be easily verified that
$\underset{k\to\infty}\lim\Theta(f_k)=\frac{1}{2}(z(x_A^{})-z(x_B^{}))$.
It only remains to prove that
$\Theta(f)\leq\frac{1}{2}(z(x_A^{})-z(x_B^{}))$ for no
homeomorphism $f\in H(X,Y)$. If such a homeomorphism existed, for
every $x\in X$ we would have
$|z(x)-z(f(x))|\leq\frac{1}{2}(z(x_A^{})-z(x_B^{}))$, and hence
$z(f(x_A^{}))\geq z(y_C^{})\geq z(f(x_B^{}))$. Therefore, points
$x\in X$ such that
$|z(x)-z(f(x))|>\frac{1}{2}(z(x_A^{})-z(x_B^{}))$ could be easily
found, contradicting our assumption.
\end{exa}

\begin{exa}[Hypothesis $(b)$ fails]
This example, introduced in \cite{DoFr04}, shows that there does
not always exist an optimal homeomorphism between two size pairs
satisfying hypotheses $(a)$ and $(c)$, but missing hypothesis
$(b)$.

\begin{figure}[htbp]
\begin{center}
\psfrag{A}{$x_A$} \psfrag{B}{$x_B$} \psfrag{fg}{$f \circ \gamma$}
\psfrag{H}{$y_D^{\,\eps}$} \psfrag{g}{$\gamma$}\psfrag{f}{$f$}
\psfrag{z}{$z$} \psfrag{X}{$X$} \psfrag{Y}{$Y$}
\includegraphics[height=7cm]{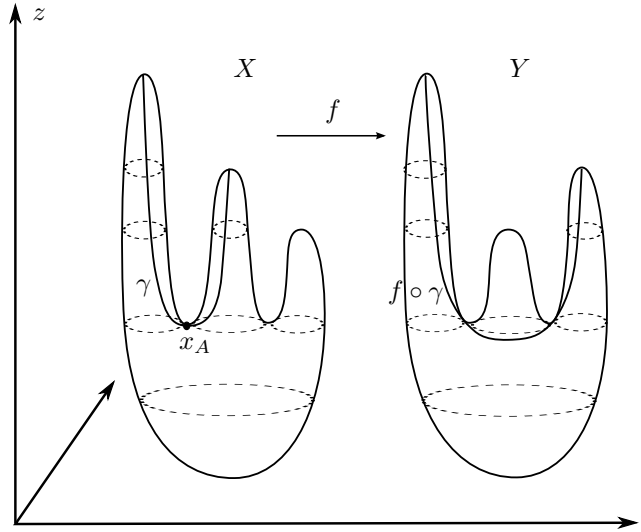}
\caption{\footnotesize{An example of two size pairs $(X, z)$ and
$(Y, z)$, whose natural pseudo-distance is zero. No optimal
homeomorphism exist between $(X, z)$, $(Y, z)$ because $X$ an $Y$
are not closed curves.}}\label{exno(c)}
\end{center}
\end{figure}

Consider the smooth surfaces $X$ and $Y$ displayed in Figure
\ref{exno(c)} and the corresponding measuring function $z$. The
dotted lines are level curves for the measuring function $z$. It
is easy to show that the natural pseudo-distance between the two
size pairs is zero. Indeed, it is possible to isotopically deform
the left surface to the right one by ``torsion'', exchanging the
positions of the two smallest humps. This deformation can be
performed by an arbitrarily small change in the values of the
height $z$. Therefore, a sequence of homeomorphisms $(f_k)$ from
$X$ to $Y$ can be constructed, such that
$\underset{k\to\infty}\lim\Theta(f_k)=0$.

However, no optimal homeomorphism exists between the two size
pairs. Suppose indeed there exists a homeomorphism $f$ such that
$\Theta(f)=0$. Consider a path $\gamma$ as in Figure
\ref{exno(c)}, chosen in such a way that, in the image of the
path, $z(x)=z(x_A)$ for no point $x \in X$ different from $x_A$.
It can be easily verified that the image of the path
$f\circ\gamma$ has to contain more than one point at which $z$
takes the value $z(x_A)$. This contradicts the assumptions, since
$\Theta(f)=0$ implies $z(f(x))=z(x)$ for every $x$ in the image of
$\gamma$.
\end{exa}

\begin{exa}[Hypothesis $(c)$ fails]
This last example shows that there does not always exist an
optimal homeomorphism between two closed curves having vanishing
natural pseudo-distance, if such curves are endowed with measuring
functions missing hypothesis $(c)$.

\begin{figure}[htbp]
\begin{center}
\psfrag{A}{$x_A$} \psfrag{B}{$x_B$}
\psfrag{C}{$y_C$}
\psfrag{D}{$x_A^{\,\eps}$}
\psfrag{E}{$x_B^{\,\eps}$} \psfrag{F}{$y_A^{\,\eps}$}
\psfrag{G}{$y_B^{\,\eps}$} \psfrag{g}{$g_{\eps}$}
\psfrag{z}{$z$} \psfrag{e}{$\eps$} \psfrag{X}{$X$} \psfrag{Y}{$Y$}
\includegraphics[height=6cm]{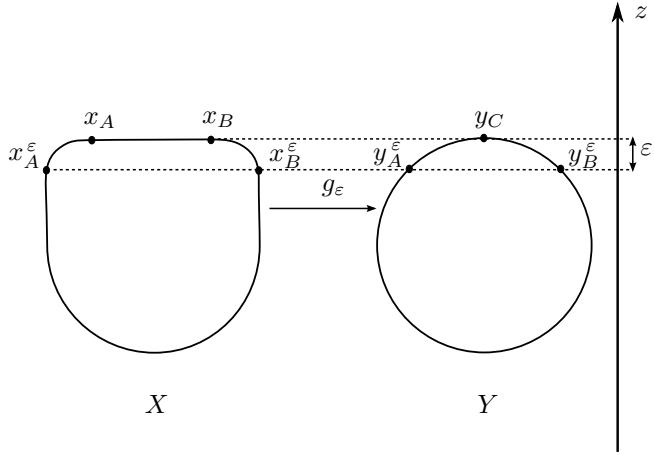}
\caption{\footnotesize{An example of two size pairs $(X, z)$ and
$(Y, z)$, whose natural pseudo-distance is zero. No optimal
homeomorphism exist between $(X, z)$, $(Y, z)$ because $z_{|_X}$
is not Morse.}}\label{exno(b)}
\end{center}
\end{figure}

Let us consider the two size pairs $(X, z)$ and $(Y, z)$ in Figure
\ref{exno(b)}, where $X$ and $Y$ are smooth closed curves. As can
be seen, the measuring function $z$ is not Morse on $X$.

We see that the natural pseudo-distance between $(X,z)$ and
$(Y,z)$ is vanishing, but an optimal homeomorphism does not exist.
Indeed, it is possible to give a sequence of homeomorphism
$(f_k)$, such that $\underset{k\to\infty}\lim\Theta(f_k)=0$, and
verify that $\Theta(f)>0$ for every homeomorphism $f\in H(X,Y)$.
Similarly to the previous example, for every $\varepsilon>0$ we
first construct a homeomorphism $g_{\varepsilon}:X\to Y$, moving
each point $x\in X$ less then or equal to $\varepsilon$ with
respect to the measuring function $z$. This can be done by
considering a homeomorphism $g_{\varepsilon}$ taking the arc
$x_{A}^{\,\varepsilon}x_{B}^{\,\varepsilon}$, containing the
segment $\overline{x_{A}x_{B}}$, to the arc
$y_{A}^{\,\varepsilon}y_C^{}y_{B}^{\,\varepsilon}$. Observe that
$z(x)=z(y_C)$ for every $x\in\overline{x_{A}x_{B}}$, and
$z(x_{A}^{\,\varepsilon})=z(x_{B}^{\,\varepsilon})=
z(y_{A}^{\,\varepsilon})=z(y_{B}^{\,\varepsilon})=z(y_C^{})-\varepsilon$.
Outside the arc $x_{A}^{\,\varepsilon}x_{B}^{\,\varepsilon}$ in
$X$ we define $g_{\varepsilon}$ by mapping every point $x$ to a
point $g_{\varepsilon}(x)$ satisfying
$z(x)=z(g_{\varepsilon}(x))$.  For every
$k\in\mathbb{N}\setminus\{0\}$ set $f_k=g_{\frac{1}{k}}$. It can
be easily verified that $\underset{k\to\infty}\lim\Theta(f_k)=0$.
However, an optimal homeomorphism $f:X\to Y$ does not exist.
Indeed, such a map should verify $\underset{x\in
X}\max|z(x)-z(f(x))|=0$, and therefore it should take each point
of the segment $\overline{x_{A}x_{B}}$ to the point $y_C$, against
the injectivity.
\end{exa}

\section{Main theorem}\label{optimal}
In this section we prove the main theorem of this paper which
states that an optimal homeomorphism exists between two closed
curves of class $C^2$, endowed with Morse measuring functions, and
whose natural pseudo-distance is zero (see Theorem \ref{mainthm}).
Roughly speaking, the proof involves the idea to construct such a
homeomorphism between the two curves as a continuous extension of
a uniformly continuous, bijective map existing between dense
subsets of the curves (see Proposition \ref{fonDense} and Remark
\ref{Extension}). The optimality is finally showed in Theorem
\ref{mainthm}.

Let us now introduce some notations and assumptions we shall adopt
in the rest of this section.

Let $(X,\p)$, $(Y,\s)$ be two size pairs, with $X$, $Y$ two $C^2$
closed curves, and $\p$, $\s$ Morse measuring functions, and
suppose that $\delta\left((X,\p),(Y,\s)\right)=0$.

It is not restrictive to assume that $X$ and $Y$ are metric
spaces, endowed with two metrics $d_X$ and $d_Y$, respectively.
Moreover, for the sake of simplicity, from now on we shall assume
that the considered curves are connected. However, note that this
last hypothesis can be weakened to any finite number of connected
components, without much affecting the following reasonings.

Let us now consider two parameterizations $h_{X}:S^1\to X$,
$h_{Y}:S^1\to Y$. The clockwise orientation on
$S^1\subset\R^2$, and the homeomorphisms $h_X,h_Y$ allow
us to induce an orientation on $X$ and $Y$, respectively. For
every $x,x'\in X$ (respectively $y,y'\in Y$), we shall denote by
$\!\!\!\xymatrix@C=0pt{{}\ar@{-} @/^/[rr]&xx'&}\!\!\!$ (resp.
$\!\!\!\xymatrix@C=0pt{{}\ar@{-} @/^/[rr]&yy'&}\!\!\!$) the
oriented path on $X$ (resp. $Y$), induced by $h_{X}$ (resp.
$h_{Y}$), from the point $x$ (resp. $y$) to the point $x'$ (resp.
$y'$), going clockwise along $S^1$, and including both $x$ and
$x'$ (resp. $y$ and $y'$).

Consider the sets
$X_{\Q}=\{x=h_{X}((\cos\theta,\sin\theta)):\theta\in\mathbb{Q}\}$
and
$Y_{\Q}=\{y=h_{Y}((\cos\theta,\sin\theta)):\theta\in\mathbb{Q}\}$,
and a sequence $\left(f_k\right)$ of $0$-approximating
homeomorphisms from $(X,\p)$ to $(Y,\s)$, i.e. such that
$\underset{k\to\infty}\lim \Theta(f_k)=0$. By using the Cantor's
diagonalization argument, and from the compactness of $X$ and $Y$,
we can assume (possibly by considering a subsequence) that there
exist $\underset{k\to\infty}\lim f_k(x)$ for every $x\in X_{\Q}$,
and $\underset{k\to\infty}\lim f^{-1}_k(y)$ for every $y\in
Y_{\Q}$. We shall set $\underset{k\to\infty}\lim f_k(x)=y_x\in Y$
for every $x\in X_{\Q}$, and $\underset{k\to\infty}\lim
f^{-1}_k(y)=x_y\in X$ for every $y\in Y_{\Q}$.

Furthermore, since homeomorphisms between closed curves can be
orientation-preserving or not, for the sake of simplicity we shall
assume (possibly by considering a subsequence) that the
orientation is maintained by each $f_k$. Indeed, if this is not
the case, we can consider a new parametrization $\hat h_Y$ having
opposite orientation with respect to $h_Y$.

Let us now set $\widetilde Y=\{y_x: x\in X_{\Q}\}$ and $\widetilde
X=\{x_{y}: y\in Y_{\Q}\}$ and denote by $X_*$, $Y_*$ the sets
$X_{\Q}\cup\widetilde X$ and $Y_{\Q}\cup\widetilde Y$,
respectively. We can define a relation $\rho\subseteq X_*\times
Y_*$ by setting
\begin{eqnarray}\label{relation}
(x,y)\in\rho\Leftrightarrow\textrm{($x\in X_{\Q}$ and
$y=y_x\in\widetilde{Y}$) or ($y\in Y_{\Q}$ and
$x=x_y\in\widetilde{X}$).}
\end{eqnarray}
\begin{rem}\label{0approx}
Note that the equality $\p(x)=\s(y)$ holds for every
$(x,y)\in\rho$. Indeed, since $\left(f_k\right)$ is a
$0$-approximating sequence, if $x \in X_{\Q}$ and
$y=y_x\in\widetilde{Y}$, then $|\p(x)- \s(y_x)|=|\p(x)-
\s(\underset{k\to\infty}\lim f_k(x))| = \underset{k\to\infty}\lim
|\p(x)-\s(f_k(x))|=0$; if $y\in Y_{\Q}$ and
$x=x_y\in\widetilde{X}$, by \emph{Remark \ref{RemDesApp}} it
follows that $|\p(x_y)- \s(y)| = |\p(\underset{k\to\infty}\lim
f_k^{-1}(y))- \s(y)| = \underset{k\to\infty}\lim |\p(f_k^{-1}(y))-
\s(y)|=0$.
\end{rem}

Following the rough outline exposed at the beginning of this
section, we shall first construct a suitable function between
$X_*$ and $Y_*$, proving that it is bijective and uniformly
continuous. Such a function can be obtained directly from the
relation $\rho$, by virtue of the following technical lemma.

\begin{lem}\label{LemmaFond}
The following statements hold:
\begin{enumerate}
    \item[$(i)$]For every real number $\varepsilon>0$, there exists
    $\eta>0$ such that, for every $(x,y),(x',y')\in\rho$ with
    $d_X\left(x,x'\right)<\eta$ (respectively
    $d_Y\left(y,y'\right)<\eta$), the inequality
    $d_{Y}\left(y,y'\right)<\varepsilon$ (resp.
    $d_X\left(x,x'\right)<\varepsilon$) holds.
    \item[$(ii)$] For every
    $x\in X_*$ (respectively $y\in Y_*$), there exists $y\in Y_*$
    (resp. $x\in X_*$) such that $(x,y)\in\rho$;
    \item[$(iii)$] For
    every $(x,y),(x',y')\in\rho$, $x=x'$ if and only if $y=y'$.
  \end{enumerate}
\end{lem}
\begin{proof}
Let us start by proving assertion $(i)$. We shall confine
ourselves to prove that for every real number $\varepsilon>0$,
there exists $\eta>0$ such that, for every $(x,y),(x',y')\in\rho$
with $d_X\left(x,x'\right)<\eta$, the inequality
$d_{Y}\left(y,y'\right)<\varepsilon$ holds. Indeed, the proof of
the other case is analogous.

We shall prove the statement by contradiction, i.e. by supposing
the existence of a real number $\bar\varepsilon>0$ such that, for
every $\eta>0$, two pairs
$(x_{\eta},y_{\eta}),(x'_{\eta},y'_{\eta})\in\rho$ exist with
$d_{X}\left(x_{\eta},x'_{\eta}\right)<\eta$ and
$d_{Y}\left(y_{\eta},y'_{\eta}\right)\geq\bar\varepsilon$. Let us
consider two sequences $((x_n,y_n))$, $((x'_n,y'_n))$ of elements in $\rho$, with
$d_X(x_n,x'_n)<\frac{1}{n}$ and $d_Y(y_n,y'_n)\geq\bar\varepsilon$
for every $n\in\mathbb{N}$.

Since $X_*=X_{\Q}\cup\widetilde{X}$, it can be assumed (possibly
by considering two subsequences) that $x_n, x'_n\in
X_{\mathbb{Q}}$ for every index $n$. Indeed, if this is not the
case, we can alternatively assume that $x_n, x'_n\in\widetilde{X}$
for every index $n$, or that $x_n \textrm{ (respectively
$x'_n$)}\in X_{\mathbb{Q}}$ and $x'_n \textrm{ (resp.
$x_n$)}\in\widetilde{X}$ for every index $n$, without much
affecting the following reasonings. Observe that our assumption
implies that $y_n, y'_n\in\widetilde{Y}$ for every index $n$.

By the compactness of $X$, we can hypothesize (possibly by
extracting a subsequence) that the sequence $(x_n)$ converges to a
point $\bar x\in X$. Obviously, it holds that
$\underset{n\to\infty}\lim d_X(x_n,x'_n)=0$ and hence $d_X(\bar
x,x'_n)\to 0$ for $n\to\infty$, that is also $(x'_n)$ converges to
$\bar x$.

Let us now consider the sequences $(y_n)$, $(y'_n)$ in
$\widetilde{Y}$. By the compactness of $Y$, we can assume that
they converge to $\bar y, \bar y'\in Y$, respectively. Moreover,
by the hypothesis $d_Y(y_n,y'_n)\geq\bar\varepsilon$ for every
index $n$, it means that $\bar y\neq\bar y'$. On the other hand,
by the continuity of $\p$ and $\s$, we have
$\underset{n\to\infty}\lim \p(x_n)=\p(\bar x)$,
$\underset{n\to\infty}\lim \p(x'_n)=\p(\bar x)$,
$\underset{n\to\infty}\lim \s(y_n)=\s(\bar y)$,
$\underset{n\to\infty}\lim \s(y'_n)=\s(\bar y')$, and by Remark
\ref{0approx} we can write $\p(\bar x)=\s(\bar y)=\s(\bar y')=c
\in\R$. In other words, we have $\bar y\neq\bar y'$ with $\s(\bar
y)=\s(\bar y')$. Since $\s$ is a Morse function, it is necessarily
non-constant on the path $\!\!\!\xymatrix@C=0pt{{}\ar@{-}
@/^/[rr]&\bar y\bar y'&}\!\!\!$, therefore there exists a point
$y'' \in \!\!\!\xymatrix@C=0pt{{}\ar@{-} @/^/[rr]&\bar y\bar
y'&}\!\!\!$ verifying $|\s\left(y''\right)-c|=C>0$ (obviously,
$y''\neq \bar y, \bar y'$). Furthermore, by recalling that
$(y_n)=(\underset{k\to\infty}\lim f_k(x_n))$ converges to $\bar y$
and $(y'_n)=(\underset{k\to\infty}\lim f_k(x'_n))$ converges to
$\bar y'$, we can find an index $N$ such that, for every $n>N$, an
index $K=K(n)$ exists with $y''\in \!\!\!\xymatrix@C=0pt{{}\ar@{-}
@/^0.8pc/[rr]&f_{k}(x_n)f_{k}(x'_n)&}\!\!\!$ for every $k>K$,
implying that $f_{k}^{-1}(y'')\in\!\!\!\xymatrix@C=0pt{{}\ar@{-}
@/^/[rr]&x_n x'_n&}\!\!\!$ for large enough indices $n$ and $k$.

Since $\underset{n\to\infty}\lim x_n=\underset{n\to\infty}\lim
x'_n=\bar x$, it can be assumed that $f_k^{-1}(y'')$ converges to
$\bar x$. Indeed, if $f_k^{-1}(y'')$ did not converge to $\bar x$,
then we could consider another point $y''' \in
\!\!\!\xymatrix@C=0pt{{}\ar@{-} @/^/[rr]&\bar y'\bar y&}\!\!\!$
(i.e. the clockwise oriented path from $\bar y'$ to $\bar y$),
verifying $|\s\left(y'''\right)-c|=C'>0$ (obviously, $y'''\neq
\bar y, \bar y'$), and such that $f_k^{-1}(y''')$ converges to
$\bar x$.

Therefore, by Definition \ref{DeSeApp} and Remark \ref{RemDesApp},
both of them in the case $d=0$, we have
$0=\underset{k\to\infty}\lim|\p(f_k^{-1}(y''))-\s(y'')|=|\p(\bar
x)-\s(y'')|=|c- \s(y'')|$, i.e. $C=0$, thus getting a
contradiction. This concludes the proof of $(i)$.

The proof of statement $(ii)$ is trivial, and directly follows by
the definition of the relation $\rho$ in (\ref{relation}).

Let us now prove $(iii)$. Let $(x,y),(x',y')\in\rho$, with $x=x'$.
This means that $d_X(x,x')<\eta$ for every real value $\eta>0$.
Since $(x,y),(x',y')\in\rho$, assertion $(i)$ implies that
$d_Y(y,y')<\varepsilon$ for every real value $\varepsilon>0$, i.e.
$y=y'$. Conversely, a similar proof can be given, by exchanging the
roles of $x,x'$ and $y,y'$ and applying once more assertion $(i)$.
\end{proof}

In plain words, assertions $(ii)$ and $(iii)$ of Lemma
\ref{LemmaFond} tell us that it is possible to define a bijective
function $\f:X_*\to Y_*$ directly from the relation $\rho$, by
setting
$$
\f(x)=y\Leftrightarrow(x,y)\in\rho.
$$
Moreover, statement $(i)$ of Lemma \ref{LemmaFond} implies that
$\f$ is uniformly continuous together with its inverse. Finally,
we observe that Remark \ref{0approx} implies $\p(x)=\s(\f(x))$ for
every $x\in X_*$. In other words, we have just proved the
following proposition.

\begin{prop}\label{fonDense}
$f_*:X_*\to Y_*$ is a bijective, uniformly continuous map with
uniformly continuous inverse, and such that $\p(x)=\s(\f(x))$ for
every $x\in X_*$.
\end{prop}

\begin{rem}\label{Extension}
Proposition \ref{fonDense} allows us to easily obtain a continuous
map from $X$ to $Y$. Indeed, it is well known that a uniformly
continuous function between two metric spaces can be {\em
univocally} extended to a continuous map between two given
completions of the spaces themselves.

In our context, $X$ and $Y$ are compact metric spaces, and hence
complete. This means that every Cauchy sequence of elements in $X$
(respectively $Y$) has a limit in $X$ (resp. $Y$). Moreover, $X$
and $Y$ are completions of their dense subsets $X_*$ and $Y_*$,
respectively. Finally, we observe that a uniformly continuous
function takes Cauchy sequences into Cauchy sequences. Following
these considerations, it is easy to show that $f:X\to Y$ with
$f_{|_{X_*}}=\f$ and $f(\underset{n\to\infty}\lim
x_n)=\underset{n\to\infty}\lim\f(x_n)$ for every Cauchy sequence
$(x_n)$ of elements in $X_*$, is a well-defined function,
continuously extending the map $\f:X_*\to Y_*$.
\end{rem}

We are now ready to give the main result of this paper.

\begin{thm}\label{mainthm}
Let $(X,\p),\ (Y,\s)$ be two size pairs, with $X,Y$ closed curves
of class $C^2$, and $\p:X\rightarrow \R$, $\s:Y\rightarrow \R$
Morse measuring functions. If
$\delta\left((X,\p),(Y,\s)\right)=0$, then there exists an optimal
homeomorphism $f:X\to Y$.
\end{thm}
\begin{proof}
We shall prove the statement by showing that the continuous
function $f$ defined in Remark \ref{Extension} is an optimal
homeomorphism from $X$ to $Y$.

Let us start by proving that $f$ is a homeomorphism. To do so, we
only need to show that $f$ is injective, since every continuous
injection between compact Hausdorff spaces is a homeomorphism
\cite[Thm. 2-103]{HoYo61}.

Let $x,x'\in X$ and suppose $x\neq x'$. Then
$x=\underset{n\to\infty}\lim x_n$, $x'=\underset{n\to\infty}\lim
x'_n$ for two suitable sequences $(x_n)$ and $(x'_n)$ of elements
in $X_*$, with $\underset{n\to\infty}\lim
x_n\neq\underset{n\to\infty}\lim x'_n$. This means that we can
choose a real number $\bar\varepsilon>0$ such that, for
sufficiently large indices $n$, we have
$d_X(x_n,x_n')>\bar\varepsilon$, allowing us to claim that
$\underset{n\to\infty}\lim\f(x_n)\neq\underset{n\to\infty}\lim\f(x'_n)$.
Indeed, the equality
$\underset{n\to\infty}\lim\f(x_n)=\underset{n\to\infty}\lim\f(x'_n)$
would imply the existence of a real number, that is
$\bar\varepsilon$, such that, for every real value $\eta>0$, two
points $\f(x_n),\f(x'_n)\in Y_*$ would exist, with
$d_Y(\f(x_n),\f(x'_n))<\eta$ and $d_X(x_n,x_n')>\bar\varepsilon$,
thus contradicting the uniform continuity of the inverse of $\f$
(see Proposition \ref{fonDense}). Hence, the assumption $x\neq x'$
implies that
$\underset{n\to\infty}\lim\f(x_n)\neq\underset{n\to\infty}\lim\f(x'_n)$,
i.e. $f(x)\neq f(x')$, thus proving the injectivity of $f$.

To conclude the proof, we still need to show the optimality of
$f$, i.e. that the equality $\p(x)=\s(f(x))$ holds for every $x\in
X$. By Proposition \ref{fonDense}, this is true when $x\in X_*$.
On the other hand, if $x\in X\setminus X_*$, then there exists a
sequence $(x_n)$ in $X_*$ converging to $x$. So, by the continuity
of $\p$, $\s$ and $f$, we can write
$\p(x)=\underset{n\to\infty}{\lim}\p(x_n)$ and
$\s(f(x))=\underset{n\to\infty}{\lim}\s(f(x_n))$. Moreover, by
recalling once more that the restriction of $f$ to $X_*$ coincides
with $\f$, and $\s\circ\f$ coincides with $\p_{|_{X_*}}$, we have
$|\p(x)-\s(f(x))|=|\underset{n\to\infty}{\lim}\p(x_n)-
\underset{n\to\infty}{\lim}\s(f(x_n))|=\underset{n\to\infty}{\lim}|\p(x_n)-\s(f(x_n))|=
\underset{n\to\infty}{\lim}|\p(x_n)-\s(\f(x_n))|= 0$.
\end{proof}

\section{Conclusions and future works}
In this paper we have proved that there always exists an optimal
homeomorphism between two size pairs $(X,\p)$, $(Y,\s)$ having
vanishing natural pseudo-distance, under the assumptions that
$X,Y$ are closed curves of class $C^2$, and $\p,\s$ are Morse
measuring functions. We point out that this result is the first
available one concerning the existence of optimal homeomorphisms
between size pairs. Indeed, until now the research has been
developed mainly focusing on the relations between the natural
pseudo-distance and the critical values of the measuring
functions, as well as on the estimation of natural pseudo-distance
via lower bounds provided by size functions. Our result opens the
way to further investigations, in order to obtain a generalization
to the case of $k$-dimensional manifolds endowed with
$\R^k$-valued measuring functions, with $k>1$. In this context, an
interesting research line appears to be, for example, to consider
measuring functions having finite preimage for each point in the
range, or characterized by a behavior analogous to that of Morse
functions in the $1$-dimensional case.

\end{document}